\documentclass[twoside,11pt]{amsart}
\numberwithin{equation}{section}
\newtheorem{thm}{Theorem}[section]

\newtheorem{prop}[thm]{Proposition}
\newtheorem{cor}[thm]{Corollary}

\theoremstyle{definition}

\theoremstyle{remark}
\newtheorem*{rem}{Remark}

\newcommand{\R}{\mathbb{R}}
\newcommand{\W}{\mathcal{W}}

\newcommand{\AK}{\mathcal{AK}}
\newcommand{\ie}{i.\thinspace{}e. }
\newcommand{\eg}{e.\thinspace{}g. }
\newcommand{\g}{\widetilde{g}}
\newcommand{\n}{\nabla}

\newcommand{\wh}[1]{\widehat{#1}}
\newcommand{\I}{\bar{i}}
\newcommand{\J}{\bar{j}}

\newcommand{\thmref}[1]{Theorem~\ref{#1}}
\newcommand{\propref}[1]{Proposition~\ref{#1}}

\newcommand{\corref}[1]{Corollary~\ref{#1}}
\begin{document}

\title[Tangent Bundles with Complete Lift Metric and Almost $hcHN$-Structure]
{Tangent Bundles with Complete Lift of the Base Metric and Almost
Hypercomplex Hermitian-Norden Structure}

\author{Mancho Manev}



\thanks{This paper is partially supported by project NI13-FMI-002
of the Scientific Research Fund, Plovdiv University, Bulgaria and
the German Academic Exchange Service (DAAD)}

\dedicatory{Dedicated to the 75th anniversary of Prof. Kostadin
GRIBACHEV}

\maketitle

{\small
\textbf{Abstract}

The tangent bundle of an almost Norden manifold and the complete
lift of the Norden metric is considered as a $4n$-manifold. It is
equipped with an almost hypercomplex Hermitian-Norden structure.
It is characterized geometrically. The case when the base manifold
is an h-sphere is considered.


}

\section{Introduction}

The investigating of the tangent bundle $TM$ of a manifold $M$
help us to study the manifold $M$. Moreover, $TM$ has own
structure closely related to the structure of $M$, which implies
mutually related geometric properties.

The geometry of the almost hypercomplex manifolds with Hermitian
metric is known (\eg~\cite{AlMa}). A parallel direction including
indefinite metrics is the developing of the geometry of the almost
hypercomplex manifolds with Hermitian-Norden metric structure. It
has a natural origination from the geometry of the $n$-dimensional
quaternionic Euclidean space.

The beginning was put by our joint works with K.~Gribachev and
S.~Dimiev in \cite{GriManDim12} and \cite{GrMa24}. More precisely
we have combined the Hermitian metric with the Norden metric with
respect to the almost complex structures of a hypercomplex
structure.

The aim of the present work is consideration of an almost Norden
manifold as a base manifold and generation of its tangent bundle
with a metric, which is a prolongation of the base metric by its
complete lift. In that way we get the tangent bundle with almost
hypercomplex Hermitian-Norden structure and characterize it.


\section{Differentiable manifolds with almost complex structures}\label{sec_ac}

\subsection{Almost complex manifolds with Hermitian metric or Norden metric} %
The notion of the \emph{almost complex manifold}
$(M^{2n},J)$ is well-known. There exists a possibility it to be
equipped with two different kinds of metrics. When $J$ acts as an
isometry on each tangent space then the manifold is an
\emph{almost Hermitian manifold}. But, in the case when $J$ acts
as an anti-isometry on each tangent space, the notion of the
so-called \emph{almost Norden manifold} (or \emph{almost
anti-Hermitian manifold}) is available. Let us consider more
precisely the latter one.

Every $n$-dimensional complex Riemannian manifold induces a real
$2n$-dimensional manifold $(M^{2n},\allowbreak{}J,g,\widetilde g)$
with a complex structure $J$, a metric $g$ and an associated
metric $\widetilde g = g (\cdot,J\cdot)$. Both metrics are
indefinite of signature $(n,n)$. This manifold is called an
\emph{almost Norden manifold} because it is introduced by
A.P.~Norden in~\cite{No}. An almost Norden manifold is of
\emph{K\"ahler type} (sometimes it is called briefly a
\emph{K\"ahler-Norden manifold}) if $J$ is parallel with respect
to the Levi-Civita connection $\nabla$ of the metric $g$. The
class of these manifolds is contained in every other class of the
almost Norden manifolds. A classification with respect to $\nabla
J$ consisting of three basic classes is given in~\cite{GaBo},
whereas the corresponding classification of almost Hermitian
manifolds is known from \cite{GrHe} and it contains four basic
classes.

\subsection{Almost hypercomplex manifolds with Hermitian-Norden structure}

Let us recall the notion of the almost hypercomplex structure $H$
on a manifold $M^{4n}$. It is the triple $H=(J_\alpha)$
$(\alpha=1,2,3)$ of anticommuting almost complex structures
satisfying the property $J_3=J_1\circ J_2$ (\cite{AlMa},
\cite{So}). Further, the index $\alpha$ runs over the range
$\{1,2,3\}$ unless otherwise stated.

A pseudo-Riemannian metric $g$ of signature $(2n,2n)$ (or a
neutral semi-Riemannian metric) on $(M^{4n},H)$ is introduced as
follows (\cite{GriManDim12})
\begin{equation}\label{1}
g(\cdot,\cdot)=g(J_1\cdot,J_1\cdot)=-g(J_2\cdot,J_2\cdot)=-g(J_3\cdot,J_3\cdot).
\end{equation}
We have called such metric a \emph{Hermitian-Norden metric}. It
generates a K\"ahler 2-form $\Phi$ and two Hermitian-Norden
metrics $g_2$ and $g_3$ by the following way
\begin{equation}\label{G}
\Phi:=g(J_1\cdot,\cdot),\qquad g_2:=g(J_2\cdot,\cdot),\qquad
g_3:=g(J_3\cdot,\cdot).
\end{equation}
Let us note that $g$ ($g_2$, $g_3$, respectively) has a Hermitian
compatibility with respect to $J_1$ ($J_3$, $J_2$, respectively)
and a Norden compatibility with respect to $J_2$ and $J_3$ ($J_1$
and $J_2$, $J_1$ and $J_3$, respectively).
On the other hand, a quaternionic inner product $<\cdot,\cdot>$
generates in a natural way the
bilinear forms $g$, $\Phi$, $g_2$ and $g_3$ by the following
decomposition: $<\cdot,\cdot>=-g+i\Phi+jg_2+kg_3$.

We have called the structure $(H,G)=(J_1,J_2,J_3;g,\Phi,g_2,g_3)$
on $M^{4n}$ an \emph{almost hy\-per\-complex Hermit\-ian-Norden
structure} or shortly an \emph{almost $hcHN$-structure}. We have
called the manifold $(M,H,G)$ an \emph{almost hypercomplex
Hermitian-Norden manifold} or shortly an \emph{almost
$hcHN$-manifold}.

It is well known, that the almost hypercomplex structure
$H=(J_\alpha)$  is a \emph{hypercomplex structure} if the
Nijenhuis tensors
$
N_\alpha(\cdot,\cdot)=
    \left[\cdot,\cdot \right]
    +J_\alpha\left[\cdot,J_\alpha \cdot \right]
    +J_\alpha\left[J_\alpha \cdot,\cdot \right]
    -\left[J_\alpha \cdot,J_\alpha \cdot \right]
$ 
vanish for each $\alpha$. Moreover, a structure $H$ is
hypercomplex if and only if two of $N_\alpha$ vanish.

We have introduced three \emph{structure $(0,3)$-tensors} of the
almost $hcHN$-manifold by \( F_\alpha (x,y,z)=g\bigl( \left(
\nabla_x J_\alpha \right)y,z\bigr)=\bigl(\nabla_x g_\alpha\bigr)
\left( y,z \right), \) where $\nabla$ is the Levi-Civita
connection generated by $g$ and $x, y, z \in T_pM$ at any $p\in
M$. Relations between the tensors $F_\alpha$ are valid, \eg
$F_1(\cdot,\cdot,\cdot)=F_2(\cdot,J_3\cdot,\cdot)+F_3(\cdot,\cdot,J_2\cdot)$
\cite{GriManDim12}.

\section{Tangent bundle with almost $hcHN$-
structure}\label{sec:(TM,H,G)}

Our purpose is a determination of an almost $hcHN$-structure
$(H,G)$ on $TM$ when the base manifold $M$ has an almost Norden
structure $(J,g,\widetilde{g})$.

We will use the horizontal and vertical lifts of the vector fields on $M$ to
get the corresponding components of the considered tensor fields on $TM$. These
components are sufficient to describe the characteristic tensor fields on $TM$
in general.

\subsection{Almost hypercomplex structure on the tangent bundle}

As it is known \cite{YaIs}, for any affine connection in $M$, the
induced horizontal and vertical distributions in $TM$ are mutually
complementary. Then we define tensor fields $J_1$, $J_2$ and $J_3$
in $TM$ by their action over the horizontal and vertical lifts of
an arbitrary vector field in $M$:
\begin{equation}\label{H}
J_1:
\left\{
\begin{aligned}
    X^H&\rightarrow -(JX)^H\\
    X^V&\rightarrow (JX)^V
\end{aligned},
\right.\qquad J_2: \left\{
\begin{aligned}
    X^H&\rightarrow X^V\\
    X^V&\rightarrow -X^H
\end{aligned},
\right.\qquad J_3: \left\{
\begin{aligned}
    X^H&\rightarrow (JX)^V\\
    X^V&\rightarrow (JX)^H
\end{aligned},
\right.
\end{equation}
where $J$ is the given almost complex structure on $M$.

By direct computations we get the following

\begin{prop}
There exists an almost hypercomplex structure $H$, defined
by~\eqref{H} in $TM$ over an almost complex manifold $(M,J)$ with
an affine connection $\nabla$. The constructed $4n$-dimensional
manifold is an almost hypercomplex manifold $(TM,H)$.
\end{prop}

Let $N_\alpha$ denotes the Nijenhuis tensor of $J_\alpha$ for each
$\alpha$ and $\wh{X},\wh{Y}\in \mathcal{T}^0_1(TM)$, \ie
\begin{equation}\label{Na-def}
N_\alpha(\wh{X},\wh{Y})=[\wh{X},\wh{Y}]
+J_\alpha[J_\alpha\wh{X},\wh{Y}] +J_\alpha[\wh{X},J_\alpha\wh{Y}]
-[J_\alpha\wh{X},J_\alpha\wh{Y}].
\end{equation}

Further, we denote the horizontal and vertical lifts by
$(\cdot)^H, (\cdot)^V \in \mathcal{T}^0_1(TM)$ of any $X, Y, Z, W
\in \mathcal{T}^0_1(M)$ at $u\in T_p M$.

For the Levi-Civita connection $\nabla$ of $M$ and its curvature
tensor $R$, then at $u\in T_p M$ we have (see also \cite{YaIs})
\begin{equation}\label{l[]}
    \begin{array}{ll}
    [X^H,Y^H]=[X,Y]^H-\{R(X,Y)u\}^V,\quad
    &[X^H,Y^V]=\left(\nabla_XY\right)^V,\\[0pt]
    [X^V,Y^V]=0,
        &[X^V,Y^H]=-\left(\nabla_YX\right)^V.
   \end{array}
\end{equation}

Using  \eqref{H}, \eqref{Na-def} and \eqref{l[]}, we get
\begin{prop}\label{Na}
Let $(M,J)$ be an almost complex manifold. Then the Nijenhuis tensors of the
structure $H$ in $TM$ for the horizontal and vertical lifts have the form
\begin{equation*}\label{N_1}
\begin{array}{l}
N_1(X^H,Y^H)=\left(N(X,Y)\right)^H\\[0pt]
\phantom{N_1(X^H,Y^H)=}+\bigl(R(JX,JY)u+JR(JX,Y)u+JR(X,JY)u-R(X,Y)u\bigr)^V,\\[0pt]
N_1(X^H,Y^V)=\bigl(\left(\nabla_{JX}J\right)(Y)-\left(\nabla_{X}J\right)(JY)\bigr)^V,\\[0pt]
N_1(X^V,Y^H)=\bigl(\left(\nabla_{Y}J\right)(JX)-\left(\nabla_{JY}J\right)(X)\bigr)^V,\qquad
N_1(X^V,Y^V)=0;
\end{array}
\end{equation*}
\begin{equation*}\label{N_2}
\begin{array}{l}
N_2(X^H,Y^H)=-N_2(X^V,Y^V)=-\bigl(R(X,Y)u\bigr)^V,\\[0pt]
N_2(X^H,Y^V)=N_2(X^V,Y^H)=-\bigl(R(X,Y)u\bigr)^H;\\[0pt]
\end{array}
\end{equation*}
\begin{equation*}\label{N_3}
\begin{array}{l}
N_3(X^H,Y^H)=\bigl(J\left(\nabla_{X}J\right)(Y)
-J\left(\nabla_{Y}J\right)(X)\bigr)^H-\bigl(R(X,Y)u\bigr)^V, \\[0pt]
N_3(X^H,Y^V)=\bigl(J\left(\nabla_{X}J\right)(Y)
+\left(\nabla_{JY}J\right)(X)\bigr)^V-\bigl(JR(X,JY)u\bigr)^H, \\[0pt]
\end{array}
\end{equation*}
\begin{equation*}
\begin{array}{l}
N_3(X^V,Y^H)=-\bigl(\left(\nabla_{JX}J\right)(Y)
+J\left(\nabla_{Y}J\right)(X)\bigr)^V-\bigl(JR(JX,Y)u\bigr)^H, \\[0pt]
N_3(X^V,Y^V)=-\bigl(\left(\nabla_{JX}J\right)(Y)
-\left(\nabla_{JY}J\right)(X)\bigr)^H+\bigl(R(JX,JY)u\bigr)^V. \\[0pt]
\end{array}
\end{equation*}
\end{prop}

The last equalities for $N_\alpha$ imply the following necessary
and sufficient conditions for integrability of $J_\alpha$ and $H$.
\begin{thm}\label{thm:Na=0}
Let $TM$ be the tangent bundle manifold
with an almost hypercomplex structure $H=(J_1,J_2,J_3)$ defined as
in~\eqref{H} and $M$ be its base manifold with almost complex
structure $J$. Then the following inter\-connections hold:
\begin{enumerate}
    \item[(i)] $(TM,J_\alpha)$ for $\alpha=1$ or $3$ is
complex if and only if $M$ is flat and $J$ is parallel; 
    \item[(ii)] $(TM,J_2)$  is complex if and only if $M$ is flat;
    \item[(iii)] $(TM,H)$ is hypercomplex if and only if $M$ is flat and $J$ is parallel.
\end{enumerate}
\end{thm}
\begin{cor}
\begin{enumerate}
    \item[(i)] $(TM,J_1)$ is complex if and only if $(TM,J_3)$ is complex.
    \item[(ii)] If $(TM,J_1)$ or $(TM,J_3)$ is complex then $(TM,H)$ is hypercomplex.
\end{enumerate}
\end{cor}

\subsection{Complete lift of the base metric on the tangent bundle}

Let us introduce a metric $\wh{g}$ on $TM$, which is the complete
lift $g^C$ on $TM$ of the base metric $g$,
 by
\begin{equation}\label{5g}
\wh{g}(X^H,Y^H)=\wh{g}(X^V,Y^V)=0,\quad
\wh{g}(X^H,Y^V)=\wh{g}(X^V,Y^H)=g(X,Y).
\end{equation}

It is known that $g^C$, associated with a (semi-)Rie\-mannian
metric $g$, is a semi-Rie\-mann\-ian metric on $TM$ of signature
$(m,m)$, where $m=\dim M$, which coincides with the horizontal
lift of $g$ when this is considered with respect to the
Levi-Civita con\-nec\-tion of $g$. This metric is introduced by
Yano and Kobayashi as $(TM,g^C)$ has zero scalar curvature and it
is an Einstein space if and only if $M$ is Ricci-flat
(\cite{YaIs}).

As it is known from \cite{YaKo}, if $\nabla$ is the Levi-Civita
connection of $M$ with respect to the pseudo-Riemannian metric
$g$, then $\nabla^C$ is the Levi-Civita connection of $TM$ with
respect to $g^C$. Since $\wh{\nabla}$ is the Levi-Civita
connection of $\wh{g}$ on $TM$ as $\nabla$ of $g$ on $M$,
then using the Koszul formula 
we obtain the
covariant derivatives of the horizontal and vertical lifts of
vector fields on $TM$ at $u\in T_p M$ as follows (see also
\cite{YaIs})
\begin{equation}\label{nabli}
\begin{array}{ll}
\wh{\nabla}_{X^H}Y^H=(\nabla_XY)^H+\left(R(u,X)Y\right)^V,\quad &
\wh{\nabla}_{X^H}Y^V=(\nabla_XY)^V,\\[0pt]
\wh{\nabla}_{X^V}Y^H=0, & \wh{\nabla}_{X^V}Y^V=0.
\end{array}
\end{equation}

After that, we calculate the components of the curvature tensor
$\wh{R}$ of $\wh{\nabla}$ with respect to the horizontal and
vertical lifts of the vector fields on $M$ and we obtain the
following non-zero components for the curvature tensors $R$ and
$\wh{R}$ as well as the Ricci tensors $\rho$ and $\wh{\rho}$,
corresponding to the metrics $g$ and $\wh{g}$ on $M$ and $TM$,
respectively (see also \cite{YaIs})
\begin{equation}\label{R}
\begin{array}{l}
\wh{R}(X^H,Y^H)Z^H=\left\{R(X,Y)Z\right\}^H+\left\{\left(\nabla_uR\right)(X,Y)Z\right\}^V,\\[0pt]
\wh{R}(X^H,Y^H)Z^V=\wh{R}(X^H,Y^V)Z^H
=\left\{R(X,Y)Z\right\}^V;\quad
\wh{\rho}(Y^H,Z^H)=2\rho(Y,Z).
\end{array}
\end{equation}

Hence, we get
\begin{cor}\label{flat}
\begin{enumerate}
    \item[(i)] $(TM,\wh{g})$ is flat if and only if $(M,g)$ is flat.
    \item[(ii)] $(TM,\wh{g})$ is Ricci flat if and only if $(M,g)$ is Ricci
    flat.
    \item[(iii)] $(TM,\wh{g})$ is scalar flat.
\end{enumerate}
\end{cor}

\begin{rem}
The results of \eqref{nabli}, \eqref{R} and \corref{flat} are
confirmed also by \cite{YaIs}, where $g$ is a 
Riemannian
metric.
\end{rem}





\subsection{Tangent bundle with almost $hcHN$-structure}

Suppose that $(M,J,g,\widetilde{g})$ is an almost complex manifold
with the pair of Norden metrics $\{g, \widetilde{g}\}$ and that
$(TM,H)$ is its almost hypercomplex tangent bundle with the
Hermitian-Norden metric structure
$\wh{G}=(\wh{g},\wh{\Phi},\wh{g}_2,\allowbreak\wh{g}_3)$ derived
(as in~\eqref{G}) from the metric $\wh{g}$ on $TM$ -- the complete
lift of $g$. The generated $4n$-dimensional manifold we denote by
$(TM,H,\wh{G})$.

Bearing in mind \eqref{H}, we verify immediately that $\wh{g}$
satisfies \eqref{1} and therefore it is valid the following

\begin{thm}\label{Sas}
The tangent bundle $TM$ equipped with the almost hypercomplex
structure $H$ and the metric $\wh{g}$, defined by~\eqref{H}
and~\eqref{5g}, respectively, is an almost hypercomplex
Hermitian-Norden manifold $(TM,H,\wh{G})$.
\end{thm}

In order to characterize the structure tensors $F_\alpha$ with
respect to $\wh{g}$ and $\wh{\nabla}$ at each $u\in T_pM$ on
$(TM,H,\wh{G})$, we use \eqref{nabli} and \eqref{H}, whence we
obtain the following

\begin{prop}\label{prop:Fa}
The nonzero components of $F_\alpha$ with respect to the
horizontal and vertical lifts of the vector fields depend on
structure tensor $F$ and the curvature tensor $R$ on
$(M,J,g,\widetilde{g})$ by the following way:
\begin{equation*}
\begin{array}{l}
F_1(X^H,Y^H,Z^H)=-R(u,X,JY,Z)-R(u,X,Y,JZ),\\[0pt]
F_1(X^H,Y^H,Z^V)=-F_1(X^H,Y^V,Z^H)=-F(X,Y,Z);\\[0pt]
F_2(X^H,Y^H,Z^V)=-F_2(X^H,Y^V,Z^H)=R(u,X,Y,Z);\\[0pt]
F_3(X^H,Y^H,Z^H)=F_3(X^H,Y^V,Z^V)=F(X,Y,Z),\quad\\[0pt]
F_3(X^H,Y^H,Z^V)=-R(u,X,Y,JZ),\quad
F_3(X^H,Y^V,Z^H)=R(u,X,JY,Z).
\end{array}
\end{equation*}
\end{prop}

By direct verification of the definition conditions of the classes
in the corresponding classifications in \cite{GrHe} and
\cite{GaBo}, we obtain
\begin{prop}\label{prop:cl}
\begin{enumerate}
    \item[(i)]
    The almost Hermitian manifold $(TM,J_1,\wh{g})$ belongs to the
    class
\( \bigl\{\mathcal{AH}\setminus
\left\{\mathcal{NK}\cup\mathcal{H}\right\}\bigr\}\cup \mathcal{K},
\) where $\mathcal{AH}$, $\mathcal{NK}$, $\mathcal{H}$ and
$\mathcal{K}$ are stand for the classes of almost Hermitian,
nearly K\"ahler, Hermitian and K\"ahler manifolds, respectively.
For the 4-dimensional case, the class of $(TM,J_1,\wh{g})$ is
restricted to the class $\mathcal{AK}$ of almost K\"ahler
manifolds.
    \item[(ii)]
The almost Norden manifold $(TM,J_{\alpha},\wh{g})$,
$(\alpha=2,3)$, belongs to the class \\ \(
\left(\W_1\oplus\W_2\oplus\W_3\right)\setminus
\left\{\left(\W_1\oplus\W_2\right)\cup\left(\W_1\oplus\W_3\right)\right\}\cup
\W_0. \)
\end{enumerate}
\end{prop}

The corresponding Lee forms are determined in an arbitrary basis \(\{e_i\}_{i=1}^{4n}\) by
\(\theta_\alpha(\cdot)=g^{ij}F_\alpha (e_i,e_j,\cdot)\). Hence, we
compute them with respect to an adapted frame.

Let $\{\widehat{e}_{A}\}=
\{e_{i}^H,e_{i}^V
\}$
be the
adapted frame
at each
point of $TM$ derived by the orthonormal basis
$\{e_{i}\}$ 
of signature $(n,n)$ at each point of $M$.
The indices
%
$i,j,\dots$ run over the ranges
$\{1,2,\dots,2n\}$, while the
indices $A,B,\dots$ the range $\{1,2,\dots,4n\}$. The summation
convention is used in relation to this system of indices.

For example, we compute as follows
\[
\begin{aligned}
\theta_{3}(Z^H)&=\wh{g}^{AB}F_{3}(\wh{e}_A,\wh{e}_B,Z^H)
=g^{ij}\{F_{3}(e^H_i,e^V_j,Z^H)+F_{3}(e^V_i,e^H_j,Z^H)\}\\[0pt]
&=
g^{ij}F_{3}(e^H_i,e^V_j,Z^H)
=g^{ij}R(u,e_i,Je_j,Z)=\rho^*(u,Z).
\end{aligned}
\]
Analogously, we have \(
\theta_{3}(Z^V)=g^{ij}F_{3}(e^H_i,e^V_j,Z^V)
=g^{ij}F(e_i,e_j,Z)=\theta(Z). \)
Thus, we obtain the following nonzero components of the Lee forms
$\theta_{\alpha}$
\[
\theta_1(Z^H)=\theta_3(Z^V)=\theta(Z),\quad%
\theta_2(Z^H)=-\rho(u,Z),\quad%
\theta_3(Z^H)=\rho^*(u,Z),
\]
where $\rho$ and $\rho^*$ are the Ricci tensor and
its associated Ricci tensor regarding $g$ and $J$. 
Therefore, we obtain
\begin{prop}\label{prop:titi}
\begin{enumerate}
    \item[(i)] $\theta_1=0$ if and only if $\theta=0$;
    \item[(ii)] $\theta_2=0$ if and only if ${\rho}=0$;
    \item[(iii)] $\theta_3=0$ if and only if
        $\theta=0$ and $\rho^*=0$.
\end{enumerate}
\end{prop}

\begin{rem}
Let us recall that the condition for vanishing of the
corresponding Lee form determines the class $\mathcal{SK}$ of the
semi-K\"ahler manifolds among the almost Hermitian manifolds and
the class $\W_2 \oplus \W_3$ among the almost Norden manifolds,
respectively.
\end{rem}

Bearing in mind \propref{prop:Fa} and \thmref{thm:Na=0}, it is
easy to conclude the following
\begin{prop}
\begin{enumerate}
    \item[(i)] $(TM,J_\alpha,\wh{g})$ for $\alpha=1$ or $3$ has a
    parallel complex structure $J_\alpha$ if and only if $(M,J,g,\widetilde{g})$ is flat and
        $J$ is parallel.
    \item[(ii)] $(TM,J_2,\wh{g})$ has a parallel complex structure $J_2$ if and only if
    $(M,J,g,\widetilde{g})$ is flat.
    \item[(iii)] $(TM,H,G)$ is a pseudo-hyper-K\"ahler manifold if and only if
        $(M,J,g,\widetilde{g})$ is flat and $J$ is parallel.
\end{enumerate}
\end{prop}

\begin{rem}
We say that a Hermitian-Norden manifold is a
\emph{pseudo-hyper-K\"ahler manifold}
, if $F_\alpha=0$ for each $\alpha\in\{1,2,3\}$, \ie the manifold
is of K\"ahler type with respect to each $J_\alpha$. According to
\cite{GriManDim12}, each pseudo-hyper-K\"ahler manifold is a flat
pseudo-Rie\-mann\-ian manifold of neutral signature.
\end{rem}

\begin{cor}
\begin{enumerate}
    \item[(i)] $J_1$ is parallel if and only if
        $J_3$ is parallel.
    \item[(ii)] If
        $J_1$ or $J_3$ is parallel then
        $(TM,H,\wh{G})$ is pseudo-hyper-K\"ahler.
\end{enumerate}
\end{cor}
\begin{cor}
\begin{enumerate}
    \item[(i)] The only complex manifolds $(TM,J_\alpha,\wh{g})$ for
        some $\alpha$ are the corresponding
        manifolds of K\"ahler type with respect to the same $J_\alpha$.
    \item[(ii)] The only hypercomplex manifolds $(TM,H,\wh{G})$ are the pseudo-hyper-K\"ahler
        manifolds.
\end{enumerate}
\end{cor}

\begin{cor}\label{cor:M}
\begin{enumerate}
    \item[(i)] If $(M,J,g,\widetilde{g})$ is flat, then $TM$ has parallel $J_2$.
    \item[(ii)] If $(M,J,g,\widetilde{g})$ is flat and its Lee form $\theta$ is zero,
    \ie $(M,J,g,\widetilde{g})\in\left\{\W_2 \oplus \W_3\right\}$,
    then $(TM,H,\wh{G})\in\mathcal{SK}(J_1)\cup\mathcal{W}_0(J_2)
    \cup\{\mathcal{W}_2\oplus\mathcal{W}_3\}(J_3)$;
    \item[(iii)] If $(M,J,g,\widetilde{g})$ is a K\"ahler-Norden manifold,
    then $(TM,J_1,\wh{g})\in\AK$.
\end{enumerate}
\end{cor}

\begin{rem}
By comparison with the Riemannian case, in \cite{To} it is shown
that $(TM,J_1,\wh{g})$ is almost K\"ahler (\ie symplectic) for any
Riemannian metric $g$ on the base manifold when the connection
used to define the horizontal lifts is the Levi-Civita connection.
\end{rem}


\subsection{Tangent bundle of h-sphere}\label{sec:sect.curv.}

Let $(M,J,g,\g)$ be a K\"ahler-Norden manifold, $\dim M = 2n
\allowbreak{}\geq 4$. Let $x, y, z, w$ be arbitrary vectors in
$T_pM$, $p \in M$. The curvatu\-re ten\-sor $R$ of type $(0, 4)$
defined by $R(x, y, z, w)\allowbreak =\allowbreak g(R(x, y)z, w)$
has the K\"ahler pro\-per\-ty $R(x, y, z, w) =-R(x, y, \allowbreak
Jz, Jw)$ in this case. This implies that the associated tensor
$R^*$ of type $(0, 4)$ defined by $R^*(x, y, z, w) =\allowbreak{}
R(x, y, z, Jw)$ has the property $R^* (x, y, z, w) =$
$\allowbreak{} -R^*(x, y, w, z)$ \cite{GaGrMi}. Therefore $R^*$
has the properties of a curvature tensor, \ie it is a
\emph{curvature-like tensor}. The following tensors are essential
curvature-like tensors:
\[
\begin{array}{l}
\pi_1(x, y, \allowbreak{}z, w) = g(y, z) g(x, w)- g(x, z) g(y,
w),\\
\pi_2(x, y, z, w) = \g(y, z) \g(x, w) -\g(x, z) \g(y, w),\\
\pi_3(x, y, z, w) = -g(y, z) \g(x, w) + g(x, z) \g(y, w)-\g(y, z)
g(x, w) + \g(x, z)\allowbreak{}g(y, w).
\end{array}
\]

Every non-degenerate 2-plane $\beta$ with respect to $g$ in
$T_pM$, $p \in M$, has the following two sectional curvatures
$k(\beta;p)=cR(x,y,y,x)$, $k^*(\beta;p)=cR^*(x,y,y,x)$,
where $c=\pi_1(x,y,y,x)^{-1}$ and $\{x,y\}$ is a basis of $\beta$.

A 2-plane $\beta$ is said to be \emph{holomorphic} (resp.,
\emph{totally real}) if $\beta= J\beta$ (resp., $\beta\perp
J\beta\neq \beta$) with respect to $g$ and $J$.

The orthonormal $J$-basis $\{e_i,e_{\bar{i}}\}$, where
$i\in\{1,2,\dots,n\}$, $\bar{i}=n+i$, $e_{\bar{i}}=Je_{i}$ of
$T_pM$ with respect to $g$ generates an orthogonal basis of null
vectors $\{e_i^H,e_{\bar{i}}^H,e_i^V,e_{\bar{i}}^V\}$ of $T_u(TM)$
with respect to $\wh{g}$. Then, the basis
$\{\xi_i,\xi_{\bar{i}},\eta_i,\eta_{\bar{i}}\}$, where
$\xi_i=\frac{1}{\sqrt{2}}(e_i^H+e_i^V)$,
$\xi_{\bar{i}}=\frac{1}{\sqrt{2}}(e_{\bar{i}}^H+e_{\bar{i}}^V)$,
$\eta_i=\frac{1}{\sqrt{2}}(e_i^H-e_i^V)$,
$\eta_{\bar{i}}=\frac{1}{\sqrt{2}}(e_{\bar{i}}^H-e_{\bar{i}}^V)$,
is orthonormal with respect to $\wh{g}$ and it has a signature
$(+,\cdots,+;-,\cdots,-;-,\cdots,-;+,\cdots,+)$ if the signature
of $\{e_i,e_{\bar{i}}\}$ is $(+,\cdots,+;-,\cdots,-)$. Then, using
\eqref{H}, we obtain $J_1\xi_{\I}=\eta_i$, $J_1\eta_{\I}=\xi_i$,
$J_2\eta_{i}=\xi_i$, $J_2\eta_{\I}=\xi_{\I}$,
$J_3\xi_{i}=\xi_{\I}$, $J_3\eta_{\I}=\eta_i$, \ie the basis
$\{\xi_i,\xi_{\bar{i}},\eta_i,\eta_{\bar{i}}\}$ is an adapted
$H$-basis.

Thus, the following basic 2-planes in $T_u(TM)$ are special with
respect to $H$ ($i\neq j$):
\begin{itemize}
    \item
$J_{\alpha}$-totally-real 2-planes ($\alpha=1,2,3$):
$\{\xi_i,\xi_j\}$, $\{\xi_i,\xi_{\J}\}$, $\{\xi_i,\eta_j\}$,
$\{\xi_i,\eta_{\J}\}$, $\{\xi_{\I},\xi_{\J}\}$,
$\{\xi_{\I},\eta_j\}$, $\{\xi_{\I},\eta_{\J}\}$,
$\{\eta_i,\eta_j\}$, $\{\eta_{i},\eta_{\J}\}$,
$\{\eta_{\I},\eta_{\J}\}$;
    \item
$J_1$-holomorphic and $J_{\alpha}$-totally-real 2-planes
($\alpha=2,3$):
    $\{\xi_{\I},\eta_{i}\}$, $\{\eta_{\I},\xi_{i}\}$;
    \item
$J_2$-holomorphic and $J_{\alpha}$-totally-real 2-planes
($\alpha=1,3$):
    $\{\eta_{i},\xi_{i}\}$, $\{\eta_{\I},\xi_{\I}\}$;
    \item
$J_3$-holomorphic and $J_{\alpha}$-totally-real 2-planes
($\alpha=1,2$):
    $\{\xi_{i},\xi_{\I}\}$, $\{\eta_{\I},\eta_{i}\}$.
\end{itemize}

The sectional curvatures $\wh{k}$ of these 2-planes and the
sectional curvatures $k_{ij}$, $k_{i\J}$, $k_{\I\J}$ and $k_{i\I}$
of the special basic 2-planes in $T_pM$ --  $J$-totally-real
2-planes $\{e_{i},e_{j}\}$, $\{e_{i},e_{\J}\}$,
$\{e_{\I},e_{\J}\}$ ($i\neq j$) and $J$-holomorphic 2-planes
$\{e_{i},e_{\I}\}$, respectively -- are related as follows:
\[
\begin{array}{lll}
\wh{k}(\xi_{i},\xi_{j})=\frac{1}{4}\left(\n_uk\right)_{ij}+k_{ij},
&
\wh{k}(\xi_{i},\xi_{\J})=\frac{1}{4}\left(\n_uk\right)_{i\J}+k_{i\J},
&
\wh{k}(\xi_{i},\eta_{j})=-\frac{1}{4}\left(\n_uk\right)_{ij},  \\
\wh{k}(\xi_{i},\eta_{\J})=-\frac{1}{4}\left(\n_uk\right)_{i\J}  &
\wh{k}(\xi_{\I},\xi_{\J})=\frac{1}{4}\left(\n_uk\right)_{\I\J}+k_{\I\J},
&
\wh{k}(\xi_{\I},\eta_{j})=-\frac{1}{4}\left(\n_uk\right)_{\I j},\\
\wh{k}(\xi_{\I},\eta_{\J})=-\frac{1}{4}\left(\n_uk\right)_{\I\J},
& \wh{k}(\eta_i,\eta_j)=\frac{1}{4}\left(\n_uk\right)_{ij}-k_{ij},
&
\wh{k}(\eta_{i},\eta_{\J})=\frac{1}{4}\left(\n_uk\right)_{i\J}-k_{i\J}, \\
\wh{k}(\eta_{\I},\eta_{\J})=\frac{1}{4}\left(\n_uk\right)_{\I\J}-k_{\I\J}.&&
\end{array}
\]

Therefore, we obtain
\begin{prop}
The manifold $(TM,H,\wh{G})$ for an arbitrary almost Norden
manifold $(M,J,g)$ has equal sectional curvatures of the
$J_1$-holomorphic 2-planes and vanishing sectional curvatures of
the $J_2$-holomorphic 2-planes.
\end{prop}

Identifying the point $z = (x^1, ..., x^{n+1}; y^1, ..., y^{n+1})$ in
$\R^{2n+2}$ with the position vector $Z$, we consider the holomorphic hypersurface
$S^{2n}(z_0; a, b)$ in the K\"ahler-Norden man\-i\-fold $(\R^{2n+2},J',g')$ defined by
$
g(z-z_0, z-z_0) = a$, $
\g(z-z_0, z-z_0) = b$,
where $(0, 0)\allowbreak\neq\allowbreak (a, b) \in\R^2$. $S^{2n}$
is called an \emph{h-sphere} with  center $z_0$ and parameters
$a$, $b$. Every $S^{2n}$, $n \geq 2$, has vanishing holomorphic
sectional curvatures and constant totally real sectional
curvatures $\nu=\frac{a}{a^2+b^2}$, $\nu^*=\frac{-b}{a^2+b^2}$
\cite{GaGrMi}. Then, according to \cite{BoGa}, the curvature
tensor of the h-sphere is
\begin{equation}\label{Rnu}
R = \nu (\pi_1 - \pi_2) + \nu^*\pi_3
\end{equation}
 and
therefore $\n R=0$. Moreover, we have  $\rho=2(n-1)(\nu
g-\nu^*\g)$, $\rho^*=2(n-1)(\nu^* \g+\nu g)$, $\tau=4n(n-1)\nu$,
$\tau^*=4n(n-1)\nu^*$, where $\rho^*=\rho(R^*)$,
$\tau^*=\tau(R^*)$. Because of the form of $\rho$,  $S^{2n}$ is
called \emph{almost Einstein}.

Let us consider the tangent bundle with almost $hcHN$-structure
$(TS,H,\widehat{G})$ with the h-sphere $(S,J,g)$ with parameters
$(a,b)$  as its base K\"ahler-Norden manifold. Then, bearing in
mind \propref{prop:Fa}, \propref{prop:cl} and \corref{cor:M}
(iii), we get that $(TS,J_1,\widehat{g})\in\AK$ and
$(TS,J_{\alpha},\widehat{g})$ $(\alpha=2,3)$ belongs to
$\left(\W_1\oplus\W_2\oplus\W_3\right)\setminus
\left(\W_i\oplus\W_j\right)$, where $i\neq j\in\{1,2,3\}$.
Moreover, according to \eqref{R} and \eqref{Rnu}, we obtain the
components of the curvature tensor $\wh{R}$ of
$(TS,H,\widehat{G})$. Then we have
\[
\begin{array}{l}
\wh{k}(\xi_{i},\xi_{j})=\wh{k}(\xi_{i},\xi_{\J})=\wh{k}(\xi_{\I},
\xi_{\J})=-\wh{k}(\eta_i,\eta_j)=-\wh{k}(\eta_{i},\eta_{\J})
=-\wh{k}(\eta_{\I},\eta_{\J})=\frac{a}{a^2+b^2},  \\
\wh{k}(\xi_{i},\eta_{j})=\wh{k}(\xi_{i},\eta_{\J})
=\wh{k}(\xi_{\I},\eta_{j})=\wh{k}(\xi_{\I},\eta_{\J})=0.  \\
\end{array}
\]
\begin{cor}
The manifold $(TS,H,\wh{G})$ for an arbitrary h-sphere  $(S,J,g)$
has constant sectional curvatures of $J_{\alpha}$-totally-real
2-planes and vanishing sectional curvatures of
$J_{\alpha}$-holomorphic 2-planes ($\alpha=1,2,3$).
\end{cor}


\bigskip

\begin{center}
\small{ \noindent \textsl{
Department of Algebra and Geometry\\
Faculty of Mathematics, Informatics and IT\\
Paisii Hilendarski University of Plovdiv\\
236 Bulgaria Blvd\\
4027 Plovdiv, Bulgaria\\
e-mail: mmanev@uni-plovdiv.bg}}
\end{center}

\end{document}